
\magnification=\magstep1
\input amssym.def
\documentstyle{amsppt}
\parindent=.5truein
\parskip=3pt
\hsize=6.5truein
\vsize=9.0truein
\raggedbottom
\TagsOnRight

\redefine \D        {\Cal D}
\define   \F        {\Cal F}
\define   \Hh       {\Cal H}
\define   \Rr       {\Cal R}
\define   \R        {{\Bbb R}}
\define   \Z        {{\Bbb Z}}
\define   \N        {{\Bbb N}}
\define   \G        {\Gamma}
\define   \ga       {\gamma}
\define   \ub       {\subset}
\redefine \b        {\beta}
\define   \a        {\alpha}
\define   \e        {\epsilon}
\redefine \sl       {\sum\limits}
\define   \bk       {\backslash}
\define   \pt       {\partial}
\define   \lm       {\lambda}
\define   \ts       {\thinspace}
\define   \arow{\hbox{$\vrule height6pt width.25pt$}
                 \!\!\lower2pt\hbox{$\longrightarrow$} }
\define    \qd       {{\unskip\nobreak\hfil\penalty50
                      \hskip2em\vadjust{}\nobreak\hfil
                      $\square$\parfillskip=0pt
                      \finalhyphendemerits=0\par}}

\count0=0

\vbox to 1truein{}
\centerline{MOMENTS OF DIRICHLET SPLINES AND THEIR}
\centerline{APPLICATIONS TO HYPERGEOMETRIC FUNCTIONS}
\vbox to 1truein{}
\centerline{Edward Neuman}
\smallskip
\centerline{Department of Mathematics}
\centerline{Southern Illinois University at Carbondale}
\centerline{Carbondale, IL 62901-4408 \ U.S.A.}
\centerline{E-mail:  GA3856\@ SIUCVMB.BITNET}
\bigskip
\bigskip

\centerline{${}^*$Patrick J. Van Fleet}
\smallskip
\centerline{Department of Mathematics}
\centerline{Vanderbilt University}
\centerline{Nashville, TN 37240 \ U.S.A.}
\centerline{E-mail:  Vanfleet\@ Athena.Cas.Vanderbilt.Edu}

\vbox to 2truein{}
AMS (MSC 1991) Subject Classification: \ 41A15, 65D07,33C45

\indent
$R$-hypergeometric function, moment, generating function,Bernstein polynomial,
\newline\indent
B\'ezier polynomial,confluent hypergeometric function, Lauricella polynomial,
\newline\indent
Appell and Lauricella functions.
\bigskip
${}^*$Current address: Department of Mathematics, Sam Houston State
University,
\newline\indent
Huntsville, TX 77340 \ U.S.A. , E-mail: Mth\_ pvf\@ shsu.edu.

\newpage
{\narrower{\bf Abstract.}  Dirichlet averages of multivariate
functions are employed for a derivation of basic recurrence formulas
for the moments of multivariate Dirichlet splines.  An algorithm  for
computing the moments of multivariate simplex splines is presented.
Applications to hypergeometric functions of several variables are discussed.

\baselineskip=24truept

\subheading{Introduction}

     In [9], H.B. Curry and I.J. Schoenberg have pointed out that
univariate $B$-splines can be constructed from volumes of slices of
convex polyhedra.  An extension of this idea to the case of
multivariate splines is due to C\. de Boor [1].  Since the geometric
construction is too complicated to be used in numerical computations,
some basic recurrence formulas for these functions have been found
(see [3,8,10,11,14,17,18,19,22,23]).  Multivariate $B$-splines (also
called simplex splines) have been studied extensively over the past
thirteen years by many researchers.  These functions have been found
useful for some applications of data fitting, computer aided geometric
design, and mathematical statistics.  In [29] the author addressed
some new problems where the simplex splines could play a prominent
role.

     Further generalizations of multivariate simplex splines appear
in [19] and [13].  Here and thereafter we call these functions
Dirichlet splines.  Our choice of terminology is motivated by the fact
that the distributional definition of the Dirichlet splines (see
(2.2)) involves the Dirichlet density function.  It has been
demonstrated that this class of splines is well designed for some
problems of mathematical statistics (see [13], [19], and references
therein).  Applications to the theory of multivariate convex functions
are reported in [26, 27].

     Recently we have noticed that there is a simple relationship
between \newline
univariate Dirichlet splines and some special functions such
as \newline
$R$-hypergeometric functions and confluent hypergeometric functions \newline
(see [28]).  In this paper, we present some results for
multivariate Dirichlet splines together with their applications to special
functions of several variables.

     The outline of this paper is as follows.  In Section 2 we
introduce notation and definitions which will be used throughout the
sequel.  In the third section we give a definition and basic
properties of a Dirichlet average of a multivariate function.  Also,
we give a recurrence relation for these averages.  In Section 4 we give
two results which play a crucial role in our subsequent considerations.
Moments of the class of splines under discussion together with two moment
generating functions are presented in Section 5.  Therein, we also give an
algorithm for computing the moments of multivariate simplex splines.
Applications to hypergeometric functions, including Appell's $F_4$ and
Lauricella's $F_B$, are discussed in Section 6.  In the same section, we
give a recurrence formula, two generating functions, and an inequality for
Lauricella polynomials.

\subheading{2.  Notation and Definitions}

     Let us introduce some notation and definitions which will be used
throughout the sequel.  By $x,y,\ldots,$ we denote elements of
Euclidean space $\R^s$ $(s \ge 1)$, i.e., $x = (x_1,\ldots,x_s)$.
Superscripts are used to number vectors.  The inner product (or dot
product) of $x,y \in \R^s$ is denoted by $x \cdot y = \sl^s_{k=1}
x_ky_k$.  For a given set $X \ub \R^s$ the symbols $[X]$ and
$\hbox{vol}_s(X)$ mean the convex hull of $X$ and the $s$-dimensional
Lebesgue measure, respectively.  We use standard multi-index notation,
i.e., for $\b \in \Z^s_+$, $|\b| = \b_1 + \cdots + \b_s$, $\b! = \b_1!
\cdots \b_s!$, $\a \le \b$ means $\a_1 \le \b_1,\ldots,\a_s \le \b_s$
$(\a \in \Z^s_+)$, $x^\b = x_1^{\b_1} x_2^{\b_2} \cdots x_s^{\b_s}$.
For $r \in \Z_+$, and $\b \in \Z^s_+$ with $|\b | = r$, the multinomial
coefficient $\binom {r}{\b}$ is defined in the usual way
$$\binom {r}{\b} = \frac{r!}{\b !}$$
By
$$
     E_n = \{ t = (t_1,\ldots,t_n) \in \R^n:
          t_j \ge 0, \ \hbox{ for all } j,\ \
          \sl^n_{j=1} t_j \le 1 \}
$$
we denote the standard $n$-simplex.  Let $\R_>$ represent the set of
all positive real numbers and  let $b = (b_0,\ldots,b_n) \in \R_>^{n+1}$.
The Dirichlet density function on $E_n$, denoted by $\phi_b$, is given
by
$$
     \phi_b(t) = \G (c)\prod^n_{i=0} [\G (b_i)]^{-1}
          t_i^{b_i-1},
                                                 \tag 2.1
$$
where $t \in E_n$, $t_0 = 1-t_1 - \cdots t_n$, and $c = b_0 + \cdots + b_n$.

     For $X = \{ x^0,\ldots,x^n \} \ub \R^s$ $(n \ge s \ge 1)$ with
$\hbox{vol}_s([X]) > 0$, the multivariate Dirichlet spline
$M(\cdot|b;X)$ is defined by requiring that
$$
     \int_{\R^s} f(x) M(x|b;X)dx
          = \int_{E_n} f(Xt)\phi_b(t)dt
                                                 \tag 2.2
$$
holds for all $f \in C_0(\R^s)$ -- the space of all multivariate
continuous functions on $\R^s$ with compact support (see [13, 19]).
Here $dx = dx_1 \cdots dx_s$, $dt = dt_1 \cdots dt_n$, $Xt =
\sl^n_{i=0} t_ix^i$.  When $b_0 = \cdots = b_n = 1$, $\phi_b(t) = n!$,
and the corresponding spline becomes a simplex spline.  The latter
spline will be denoted by $M(\cdot|X)$.  When $s = 1$, we will write
$m(\cdot|b;Z)$ instead of $M(\cdot|b;X)$, where $Z = \{ z_0,\ldots,z_n
\} \ub \R$.  In this case (2.2) becomes
$$
     \int_{\R} h(u)m(u|b;Z)du
          = \int_{E_n} h(Zt)\phi_b(t)dt
                                                 \tag 2.3
$$
$(h \in C_0(\R))$.

\subheading{3.  Dirichlet Averages}

     The purpose of this section is two-fold.  We give a definition of
Dirichlet averages of multivariate functions.  Next we prove a
recurrence formula for these averages.  This result has an immediate
application in Section 5.  For the reader's convenience, let us
recall a definition of the Dirichlet average of a univariate function
$h \in C_0(\R)$.  Assume that the set $Z = \{ z_0,\ldots,z_n \} \ub
\R$ is such that $\min\{z_j: 0 \le j \le n \} < \max \{z_j: 0 \le j
\le n \}$.  For $b \in \R_>^{n+1}$, the Dirichlet average of $h$,
denoted by $H(b;Z)$, is given by
$$
     H(b;Z) = \int_{E_n} h(Zt)\phi_b(t)dt
                                                 \tag 3.1
$$
(see [6]).  Comparison with (2.3) yields
$$
     H(b;Z) = \int_{\R} h(u)m(u|b;Z)du.
                                                 \tag 3.2
$$
We list below some elementary properties of $H(b;Z)$.
\item{(i)}  $H(b_0,\ldots,b_n;z_0,\ldots,z_n)$ is symmetric in indices
     $0,1,\ldots,n$ (see [6, Thm. 5.2--3]).
\item{(ii)}  A vanishing parameter $b_i$ can be omitted along with the
     corresponding variable $z_i$ (see [6, (6.3--3)]).
\item{(iii)}  Equal variables can be replaced by a single variable if
     the corresponding parameters are replaced by their sum  (see
     [6, Thm. 5.2--4]).

     We now introduce the Dirichlet average of $f \in C_0(\R^s)$.  For
$X \ub \R^s$ with $\hbox{vol}_s([X])~>~0$ and $b \in \R_>^{n+1}$ $(n
\ge s \ge 1)$, the Dirichlet average of $f$, denoted by $\F(b;X)$, is
given by
$$
     \F(b;X) = \int_{E_n} f(Xt)\phi_b(t)dt,
                                                 \tag 3.3
$$
where $Xt$ and $\phi_b$ have the same meaning as in Section 2.
Comparison with (2.2) shows that
$$
     \F(b;X) = \int_{\R^s} f(x)M(x|b;X)dx.
                                                 \tag 3.4
$$

     It is clear that the properties (i)--(iii) are also valid for the
average $\F$.  In particular, property (iii), when applied to
$M(\cdot|b;X)$ yields
$$
     \aligned
     M(\cdot|b;X) = &M(\cdot|x^0,\ldots,x^0,
                               \ldots,x^n,\ldots,x^n)\\
\vspace{-3\jot}
              &\qquad(b_0-\hbox{times})\quad\
                     (b_n-\hbox{times})
     \endaligned
                                                 \tag 3.5
$$
provided the $b$'s are positive integers.  (See also [13, 19].)  The
spline on the right hand side of (3.5) is a multivariate simplex
spline with coalescent knots (see [14] for a detailed analysis of this
class of splines).

     Before we state and prove the first result of this section, let us
introduce more notation.  By $e_j$ $(0 \le j \le n)$, we denote the
$j$th coordinate vector in $\R^{n+1}$.  For $f \in C_0^1(\R^s)$ define a
function $f_j$ as follows:
$$
     f_j(x) = D_{x^j-x} f(x)
                                                 \tag 3.6
$$
Here $D_yf$ denotes the directional derivative of $f$ in the direction
$y \in \R^s$, i.e.,
$$
     D_y f(x) = \sl^s_{k=1} y_k
          \frac{\pt}{\pt x_k} f(x).
$$
We are now ready to prove the following.

\demo{THEOREM 3.1 ([7])} Let $X = \{ x^0,\ldots,x^n \} \ub \R^s$ $(n \ge s
\ge 1)$ be such that $\hbox{vol}_s([X]) > 0$.  Further, let $f \in
C^1_0(\R^s)$ and let the vector $b \in \R_>^{n+1}$ be such that $b_j
\ge 1$ for some $0 \le j \le n$.  Then the following identity
$$
     (c-1) \F(b;X)
          = (c-1)\F(b-e_j;X) + \F_j(b;X)
                                                 \tag 3.7
$$
is valid.  Here $\F_j$ denotes the Dirichlet average of the function
$f_j$.
\enddemo

\demo{REMARK}  The proof presented below bears no resemblance to what was
done in [7,\ Thm. 3].  In this paper, the author has established (3.7) using
generalized Euler-Poisson partial differential equations.
\enddemo
\demo{Proof}  In order to establish the identity (3.7) we employ the
following one
$$
     (c-1)H(b;Z) = (c-1)H(b-e_j;Z)
          + H_j(b;Z).
                                                 \tag 3.8
$$
Here $H_j$ stands for the Dirichlet average of the function
$h_j(u) = (z_j-u)h'(u)$, $h \in C^1_0(\R^s)$.  The relation (3.8)
readily follows from (5.6--13) in [6].  Application of (3.2) to (3.8)
yields
$$
     \aligned
     (c-1) \int_{\R} h(u)m(u|b;Z)du &= (c-1) \int_{\R} h(u)
          m(u|b-e_j;Z)du      \\
     &\qquad \quad + \int_{\R} (z_j-u)h'(u) m(u|b;Z)du.
     \endaligned
                                                 \tag 3.9
$$
We will lift (3.9) to the case of multivariate functions.  To this
aim we shall employ the following formula
$$
     \int_{\R} h(u) m(u|b;Z)du = \int_{\R^s} h(\lm \cdot x)
          M(x|b;X)dx,
                                                 \tag 3.10
$$
where now $Z = \{ \lm \cdot x^0,\ldots,\lm \cdot x^n \}$, $\lm \in
\R^s \bk \{0\}$, and $h(\lm \cdot x)$ is a ridge function (or plane
wave).  Since the proof of (3.10) is similar to that presented in
[23], p\. 496, we omit further details.  Application of (3.10) to
(3.9) yields
$$
     \aligned
     (c-1)\int_{\R^s} h(\lm \cdot x)&M(x|b;X)\ts dx
          = (c-1) \int_{\R^s} h(\lm \cdot x) M(x|b-e_j;X)dx \\
     & + \int_{\R^s} (\lm \cdot x^j - \lm \cdot x)
          h'(\lm \cdot x)M(x|b;X)dx.
    \endaligned
                                                 \tag 3.11
$$
We appeal now to the denseness of ridge functions (these functions
form a dense subset in $C^1_0(\R^s)$) to conclude that the above
identity is valid for any multivariate function $f \in C^1_0(\R^s)$,
see [20].  Substituting $h(\lm \cdot x) = f(x)$ into (3.11), we
obtain the assertion and the proof is completed. \qd
\enddemo
\newpage
\demo{COROLLARY 3.2}  Along with the hypotheses of Theorem 3.1,
assume that for some $0 \le i$, $j \le n$, $1 \le k \le s$, that
$x^i_k \not= 0$, $x^j_k \not= 0$, and $b_i \ge 1$, $b_j \ge 1$.
Then
$$
     (c-1)[\F(b-e_j;X)-\F(b-e_i;X)]
         + \F_j(b;X)-\F_i(b;X) = 0
                                         \tag 3.12
$$
and
$$
     \aligned
     (c-1)(x^i_k-x^j_k) \F(b;X)&=
     (c-1)[x^i_k\F(b-e_j;X)-x^j_k\F(b-e_i;X)] \\
         &\quad+ x^i_k\F_j(b;X)-x^j_k\F_i(b;X).
     \endaligned
                                         \tag 3.13
$$
\enddemo

\demo{REMARK}  (3.12) and (3.13) are generalizations of
Exercise 5.9-6 in [6] with the latter being an extension of Zill's
identity for $R$-hypergeometric functions.

     Since the proof of (3.12) and (3.13) follows the lines introduced
in [6, p.\ 305], we omit further details.
\enddemo

\subheading{4.  Auxiliary Results}

     Our first result reads as follows.

\demo{PROPOSITION 4.1}  Let $p(x)$ be an affine function on $\R^s$.  Then
$$
     p(x)M(x|b;X) = \sl^n_{i=0}
          w_ip(x^i)M(x|b+e_i;X)
                                                 \tag 4.1
$$
provided the splines $M(x|b+e_i;X)$, $0 \le i \le n$, are continuous
at $x \in \R^s$.  Here $w_i = b_i/c$, $i=0,1,\ldots ,n$.
\enddemo

\demo{Proof}  We need the following identity for the Dirichlet density
function [6, (4.4-8)]:
$$
     t_i\phi_b(t) = w_i\phi_{b+e_i}(t)
                                                 \tag 4.2
$$
$(t \in E_n$; $0 \le i \le n)$.  Since $t_0 + \cdots + t_n = 1$, (4.2) gives
$$
     \phi_b(t) = \sl^n_{i=0} w_i\phi_{b+e_i}(t).
$$
Multiplying both sides by $f(Xt)$ and next integrating over $E_n$, we
obtain by virtue of (2.2)
$$
     \int_{\R^s} f(x)M(x|b;X)dx
          = \int_{\R^s} f(x) \sl^n_{i=0}
          w_iM(x|b+e_i;X)dx.
$$
Hence, (4.1) follows when $p(x) = 1$.  To complete the proof we
utilize (4.2) again.  Multiplying both sides by $x^i_l$ and next
summing over $i$, we obtain
$$
     [(Xt)_l] \phi_b(t) = \sl^n_{i=0}
          w_ix^i_l\phi_{b+e_i}(t).
$$
Here $(Xt)_l$ denotes the $l$th component of $Xt$.  This leads to the
following integral relation:
$$
     \int_{\R^s} f(x) x_l M(x|b;X)dx
          = \int_{\R^s} f(x) \biggl [
          \sl^n_{i=0} w_ix^i_l M(x|b+e_i;X)
          \biggr ] dx,
$$
which proves (4.1) when $p(x) = x_l$, $(1 \le l \le s)$. \qd
\enddemo

     Micchelli [22] gave a different proof of (4.1) for simplex
splines.  For this class of splines, identity (4.1) is called the
``degree elevating formula''.  A special case of (4.1) appears in [14].

     For our further aims, we recall a definition of the
$R$-hypergeometric function in the real case.  Let the set $Z=\{z_0,
\ldots ,z_n\} \ub \R^{n+1}$ be such that $0 \notin [Z]$.  Further, let
$b\in\R^{n+1}_>$.  The $R$-hypergeometric function $R_{-a}(b;Z)$,
$(a \in \R)$, is given by
$$
     R_{-a}(b;Z) = \int_{E_n} (Zt)^{-a} \phi_b(t)dt
                                                 \tag 4.3
$$
(see [6]).  When $-a \in \N$, the restriction $0\notin [Z]$ can be dropped.
Comparison with (3.1) shows that the $R_{-a}$ is the
Dirichlet average of the power function $u^{-a}$.  It is worthy to
mention that the Gauss hypergeometric function ${}_2F_1$, Lauricella's
hypergeometric function $F_\D$, the Gegenbauer polynomials, and the
elliptic integrals in the Legendre form can all be represented in
terms of the function $R_{-a}$.  Combining (4.3) and (2.3) gives
$$
     R_{-a}(b;Z) = \int_{\R} u^{-a} m(u|b;Z)du.
                                                 \tag 4.4
$$

     For future use, let us record a very useful formula for \newline
$R$-hypergeometric functions (see [6], Thm. 6.8--3)
$$
     R_{-a}(b;Z) = \prod^n_{j=0} z_j^{-b_j}
          R_{a-c}(b;Z^{-1}),
                                                 \tag 4.5
$$
where $Z^{-1}:= \{ z_0^{-1},\ldots,z_n^{-1}\}$ $(z_j > 0$, for all
$j$; $c \not= 0,-1,-2,\ldots)$.  This important result is commonly
referred to as Euler's transformation.

     We close this section with the following:

     \demo{PROPOSITION 4.2}  Let $a \in \R$ and let the vector $\lm
\in \R^s$ be such that $\lm \cdot x^j < 1$, $j=0,\ldots ,n$.  Then
$$
     \int_{\R^s} (1-\lm \cdot x)^{-a}M(x|b;X)dx
          =  R_{-a}(b;Y),
                                                 \tag 4.6
$$
where
$$
Y = 1-\lm\cdot X = \{ (1-\lm \cdot x^0),\ldots,(1-\lm \cdot x^n)\}.
                                                            \tag 4.7
$$
\enddemo

\demo{Proof} Substituting $Z = Y$ into (4.4) we obtain
$$
\aligned
   R_{-a}(b;Y) &= \int_{\R}u^{-a}m(u|b;1-\lm\cdot X)du \\
               &= \int_{\R}(1-u)^{-a}m(u|b;\lm\cdot X) \\
               &= \int_{{\R}^s}(1-\lm\cdot X)^{-a}M(x|b;X)dx.
\endaligned
$$
In the last step we have used (3.10).
\qd
\enddemo

\noindent When $a = c$, (4.6) becomes Watson's identity (see
[31])
$$
     \int_{\R^s} (1-\lm \cdot x)^{-c} M(x|b;X) dx
          = \prod^n_{j=0} (1-\lm \cdot x^j)^{-b_j}.  \tag 4.8
$$
The above identity follows by applying (4.5) to the right side of (4.6) and
using $R_0 = 1$.

\noindent
An alternative proof of (4.8) appears in [13] (see also [19]
for some comments concerning this identity).
Subject: file 3

\subheading{5.  Moments of Multivariate Dirichlet Splines}

     A motivation for the investigation of the moments of Dirichlet
splines has its origin in two mathematical disciplines.  It is well
known that the spline $M(\cdot|b;X)$ is a probability density function
on $\R^s$.  We feel that the results of this section can be applied
to some problems in mathematical statistics.  A second area of
possible applications is the theory of special functions.  We have
already mentioned that some important special functions can be
represented by the $R$-hypergeometric functions.  For particular
values of the $a$ parameter and the $b$ parameters in (4.4), this function
becomes a complete symmetric function.  For particular values of the
$z$-variables, (4.4) gives an integral formula for the $q$-binomial
coefficients (Gaussian polynomials).  (See [25] for more details).

     In this section, we derive recurrence formulas for the moments of
multivariate Dirichlet splines.  Also, we discuss implementation of
these results in the case when $b = (1,\ldots,1) \in \R^{n+1}$.  For
related results when $s = 1$, see [24].  We employ the multi-index
notation introduced in Section 2.

     For $\beta \in \R^s$, we define the moment of order
$|\b|\ $ $(|\b| = \b_1 + \cdots + \b_s)\ $ of
$M(\cdot|b;X)$ as follows:
$$
     m_\b (b;X) = \int_{\R^s} x^\b M(x|b;X) dx
                                                 \tag 5.1
$$
provided $0_s \notin [X]$, $0_s$ -- the origin in $\R^s$.
When $\b\in\Z_+^s$,  this restriction is nonessential.
In the case of the simplex spline, we shall omit the vector $b$ and
write $m_\b(X)$ instead of $m_\b(b;X)$.  Also, let $d_l$ stand for the
$l$th coordinate vector in $\R^s$.

     We are now ready to state and prove the following.

     \demo{THEOREM 5.1}  Let the weights $w_0,\ldots,w_n$ be the
same as in Proposition 4.1.  Then
$$
m_{\b} (b;X) = \sl^n_{i=0} w_im_{\b}(b+e_i;X),
                                                 \tag 5.2
$$
$$
     m_{\b+d_l} (b;X) = \sl^n_{i=0} w_ix^i_l
          m_\b(b+e_i;X),
                                                 \tag 5.3
$$
for all $l = 1,2,\ldots,s$.  Moreover, if $\hbox{vol}_s([X]) > 0$ and
$b_j \ge 1$, for some $0 \le j \le n$, then
$$
     (c+|\b|-1) m_\b(b;X)
          = (c-1)m_\b(b-e_j;X) + \sl^s_{l=1}
          \b_l x^j_l m_{\b-d_l} (b;X).
                                                 \tag 5.4
$$
If for some $0 \le i$, $j \le n$, $1 \le k \le s$,
$x^i_k \not= 0$, $x^j_k \not= 0$, and $b_i \ge 1$, $b_j \ge 1$,
then
$$
     (c-1)[m_\b(b-e_j;X)-m_\b(b-e_i;X)]
          + \sum_{k=1}^s \b_k(x^j_k-x^i_k)
          m_{\b-d_k}(b;X) = 0
                                                 \tag 5.5
$$
and
$$
     \aligned
       (c+|\b|-1) (x^i_k-x^j_k) &m_\b(b;X) =
          (c-1)[x^i_k m_\b(b-e_j;X)-x^j_k m_\b(b-e_i;X)]  \\
         & \qquad + \sum_{l=1}^s \b_l W_{k,l}
          m_{\b-d_l}(b;X),
         \endaligned
                                                 \tag 5.6
$$
where $W_{k,l} = \hbox{det }\bmatrix x^i_k &x^j_k \\
                                x^i_l &x^j_l \endbmatrix$.
\enddemo

     \demo{REMARK}  When $b_j = 1$, formula (5.4) holds true provided
$n > s$.
\enddemo

\demo{Proof}  In order to establish the recursion (5.2) and (5.3), we
substitute $p(x) = 1$ and $p(x) = x_l$ respectively, into (4.1)
and next integrate over $\R^s$.  For the proof
of (5.4), we utilize formula (3.7) with $f(x) = x^\b$.  The resulting
equation, together with (3.6) and (5.1), yields the assertion.
Formulas (5.5) and (5.6) follow immediately from Corollary 3.2 with
$f(x) = x^\b$. \qd

\enddemo

     We now give two moment generating functions.  The first generating
function involves the confluent hypergeometric function $S$.  Following
[6, (5.8-1)], we define
$$
S(b;Z) = \int_{E_n} exp(Zt)\phi_b(t)\ts dt,
                                                   \tag 5.7
$$
$b\in \R^{n+1}_>$, $Z=\{z_0,\ldots ,z_n\}$.  Use of (3.1) and (3.2) gives
$$
S(b;Z) = \int_{\R} exp(u)m(u|b;Z)\ts du.
$$
Letting $Z=\lm\cdot X = \{\lm\cdot x^0,\ldots , \lm\cdot x^n\}$,
$\lm\in\R^s\bk \{0\}$, and next using (3.10), we arrive at
$$
S(b;\lm\cdot X) = \int_{\R^s} exp(\lm\cdot x)M(x|b;X)\ts dx.
                                                               \tag 5.8
$$

To obtain the first moment generating function, we expand $exp(\lm\cdot x)$
into a power series.  Applying the multinomial theorem to powers of
$\lm\cdot x$ and next integrating the corresponding power series one term at
a time, we obtain by virtue of (5.8) and (5.1)
$$
S(b;\lm\cdot X) = \sum \frac{\lm^j}{j!}\ts m_j(b;X)
                                                               \tag 5.9
$$
where the summation extends over all multi-indices $j\in\Z^s_+$.

     It is worthy to mention that the hypergeometric function $S(b;\lm\cdot X)$
can be expressed as a divided difference of $exp(z)$ provided that \newline
$b_0, \ldots ,b_n \in \Z_+$.  We have,
$$
k!\ts S(b;\lm\cdot X) = [\lm\cdot x^0(b_0),\ldots ,\lm\cdot x^n(b_n)]e^z
                                                               \tag 5.10
$$
where $k=b_0+ \cdots b_n -1$.  Here the symbol $\lm\cdot x^i(b_i)$ means
that the knot set $\lm\cdot x^i$ is repeated $b_i$ times.  Formula (5.10)
readily follows from (5.8), (5.7), and the Hermite-Genocchi formula for
divided differences.

     A second generating function is given by
$$
R_{-a}(b;Y)  = \sum \lm^j\ts \frac{(a,|j|)}{j!}\ts m_j(b;X)  \tag 5.11
$$
($|\lm\cdot x^i| < 1$, for all $i$), where the summation extends over all
multi-indices $j\in\Z^s_+$.  Here, $a\in\R$, the set $Y$ is given in (4.7),
and $(a,l)$ stands for the Appell symbol, i.e., $(a,0)=1$, $(a,l)=
a(a+1)\cdots (a+l-1)$, $l\in\N$.  In order to establish (5.11) we expand
$(1-\lm\cdot x)^{-a}$ into a power series and next utilize the multinomial
theorem to obtain
$$
(1-\lm\cdot x)^{-a} = \sl^{\infty}_{l=0} (a,l) \sum_{|j|=l}
\frac{\lm^j}{j!}\ts x^j
                                                               \tag 5.12
$$

To complete the proof we substitute (5.12) into (4.6) and next integrate
term by term.  Applications of (5.9) and (5.11) are discussed in the next
section.

     We shall now turn our attention to the case of multivariate
simplex splines.  To this end, let $\b \in \Z^s_+$.
In the case under discussion, the formulas of (5.3) and (5.4) of
Theorem 5.1 take the form
$$
     \align
     &m_{\b+d_l}(X) = \frac{1}{n+1}
          \sl^n_{i=0} x^i_l m_\b(X^i)
                                                 \tag 5.13\\
     (n+|\b|)&m_\b(X) = n [m_\b(X_j)]
          + \sl^s_{l=1} \b_l x^j_l m_{\b-d_l}(X),
                                                 \tag 5.14
     \endalign
$$
$1 \le j \le n$; $1 \le l \le s$.  Here $X^i = X \cup \{ x^i \}$, and
$X_i = X \backslash \{ x^i \}$,
$0 \le i \le n$.  The set $X^i$ appears on the right hand side of
(5.13) because of (3.5).  A closer look at (5.14) shows that the
recursion is in two directions.  That is, given $X =
\{x^0,\ldots,x^n\} \ub \R^s$ $(n > s)$, to compute $m_\b(X)$, we need
the moment of order $|\b|$ for the knot set consisting of one less
vector than $X$, and also $s$ moments of order $|\b|-1$.

     Define a set $\overline X_k = \{ x^0,\ldots,x^k \}$, $k =
s,s+1,\ldots,n$, and note \newline
$m_\b(\overline X_k) = 1$ when $|\b| = 0$,
$k = s,s+1,\ldots,n$.  To employ (5.14), we must precompute certain moments
of the form $m_{\b} (\overline X_s)$, $|\b |>0$, and $m_{d_l}(\overline X_k)$,
$k=s,\ldots ,n$.  Let us note that
$$
m_{d_l}(\overline X_k) = \frac{1}{k+1} \sl^k_{j=0}{x_l}^j
$$
follows immediately from the defining equation (2.2).

In order to compute the moments $m_\b(\overline X_s)$, $|\beta | >0$, we first
introduce some new notation and next appeal to the proposition that follows.

     Let $t=(t_0,t_1,\ldots ,t_n) \ts\in \R^{n+1}$ with $(t_1,\ldots ,t_n)
\ts\in E_n$, and $t_0 = 1 - \sl^n_{i=1} t_i$.  Also,
let $l = (l_0,\ldots ,l_n) \ts\in\Z^{n+1}_+$, with $m=|l|$.  Then we
define the Bernstein polynomial by
$$
{B_l}^m (t) = {m \choose l} t^l.
$$

      For given coefficients $\{{p_l}\}$, $|l| = m$, we shall call any
polynomial of the form
$$
q(t) = \sum_{|l|=m} p_l {B_l}^m (t)
$$
a B\'ezier polynomial.  It is well known that any such $q$ may be
stably and efficiently evaluated using deCasteljau's algorithm
(see [2],[12]).

     Proposition 5.2 illustrates that we may indeed incorporate deCasteljau's
algorithm when evaluating $m_\b(\overline X_s)$.

\newpage
\demo{PROPOSITION 5.2}  Let $n\geq s\geq 1$, $\b\ts\in \Z^s_+$, and
$X = \{x^0,\ldots ,x^n\}\ub \R^s_+$ with
$vol_s([X]) > 0$.  Let $y^i = (x^0_i,\ldots , x^n_i)
\ts\in\R^{n+1}_+$ and set $g_i = |y^i|$, \newline $i=1,\dots ,s$.  Then
$$
m_{\b}(X) = \frac{g^{\b} n!}{(|\b|+n)!}
\sum_{|k_1|={\beta}_1} {B_{k_1}}^{\beta_1}({\tilde y}^1)
\cdots
\sum_{|k_s|={\beta}_s} {B_{k_s}}^{\beta_s} ({\tilde y}^s){\eta}!.
                                                         \tag 5.15
$$
Here,
$k_i = (k_{i_0},\ldots ,k_{in})\ts\in\Z^{n+1}_+$,
${\eta}! = {\eta_0}!\cdots{\eta_n}!$, with $\eta_j = \sl^s_{i=1} k_{ij}$,
\newline
$0\leq j\leq n$, ${\tilde y}^i = (1/g_i)y^i$, $i=1,\ldots ,s$.
\newline
$g=(g_0,\ldots ,g_n)\ts\in \R^{s+1}_+$.
\enddemo

\demo{Proof}  We use the defining relation (2.2) to write:
$$
\aligned
m_{\b}(X) &= n! \int_{E_n} \prod_{j=1}^s (t_0 x^{0}_j +
\cdots + t_n x^{n}_j )^{\b_j} \ts dt \\
& \qquad = n! \sum_{|k_1|=\b_1} \binom{\b_1}{k_1} (y^1)^{k_1}
\cdots \sum_{|k_s|=\b_s} \binom{\b_s}{k_s} (y^s)^{k_s}
\int_{E_n} t^{\eta} \ts dt,
\endaligned
$$
where $\eta$, $y^i$, $k_i$, $i=1,\ldots ,s$, are given in Proposition 5.2.
Using 4.3-4 in [6] to simplify the integral in the above identity, we have
$$m_{\b}(X) = n! \sum_{|k_1|=\b_1} \binom{\b_1}{k_1} (y^1)^{k_1}
\cdots \sum_{|k_s|=\b_s}\binom{\b_s}{k_s} (y^s)^{k_s}\frac{{\eta}!}
{(|\b|+n)!}.
                                                            \tag 5.16
$$

Now ${\tilde y}^j$, $j=1,\ldots ,s$, given in Proposition 5.2 can be viewed
as the barycentric coordinates of some point in $\R^n$.  Inserting
${\tilde y}^j$'s into (5.16) and scaling by $g^{\beta}$ gives the desired
result. \qd
\enddemo

     Thus all moments of the form $m_{\b}(\overline X_s)$ can be expressed
as a nested sum of B\'ezier polynomials and subsequently may be
evaluated using deCasteljau's algorithm.  It should be noted that while
both de Casteljau's algorithm and (5.14) are possible candidates for the
task of computing $m_{\b}(\overline X_k)$, $k=s+1, \ldots ,n$, the latter
scheme requires the evaluation of fewer terms at each recursion step and
is thus the preferred choice.  In order to summarize the procedure
for evaluating the moments of simplex splines in Algorithm 5.3, we introduce
$d_0 = (0,\ldots ,0) \in\R^s$.

\demo{ALGORITHM 5.3}

\noindent Given $X = \{ x^0,\ldots,x^n \}\ub\R^s$, $(n \ge s)$ and $\b \in
\Z^s_+$ with $|\b| > 1$, this algorithm generates the moment $m_\b(X)$
of the simplex spline $M(\cdot | X)$.

\baselineskip=16truept
\noindent
\hbox to .3in{\hfill 1.}\enspace
               $\a := d_0$

\noindent
\hbox to .3in{\hfill 2.}\enspace
        For $k = s$ to $n$  \newline
\hbox to .5truein{}  \arow\ $m_\a(\overline X_k) = 1$

\noindent
\hbox to .3in{\hfill 3.}\enspace
                For $|\a| = 1$ to $|\b|$,
                      $\a \in \Z^s_+$, $\a \le \b$
                \newline
\hbox to .4in{}  \arow\ Use (5.15) to express $m_\a(\overline X_s)$ in
                        terms of B\'ezier polynomials \newline
\hbox to .8in{}         and evaluate using DeCasteljau's algorithm.

\noindent
\hbox to .3in{\hfill 4.}\enspace
            For $k = s+1$ to $n$ \newline
\hbox to .4in{} \vrule \hbox to .2in{}
            For $\a \in \Z^s_+$, $\a \le \b$
                    \newline
\hbox to .4in{} \arow\ \arow Compute $m_\a(\overline X_k)$
                    using (5.14).
\enddemo\baselineskip=24truept

     We close this section with a remark that this algorithm is
numerically stable if $x^j > 0$ for all $j = 0,1,\ldots,n$.

\subheading{6. Applications to Hypergeometric Functions}

     In this section we demonstrate a relationship between Dirichlet
splines and an important class of hypergeometric functions of several
variables.  We will deal mainly with Appell's $F_4$ and Lauricella's $F_B$.
The link between these classes of functions is provided by another integral
average which is commonly referred to as a double Dirichlet average (see [4]
for more details).  Throughout the sequel the double Dirichlet average of a
continuous univariate function $h$ will be denoted by $\Hh$.

     Let $X\in\R^{s\times (n+1)} (n\ge s\ge 1)$.  Further, let $u=(u_1,\ldots
 ,u_s)$ be an ordered $s$-tuple of nonnegative numbers with $u_1 + \cdots +
u_s = 1$, and similarly $v=(v_0,\ldots ,v_n)$.  We define
$$
      u\cdot Xv = \sl^s_{i=1} \sl^n_{j=0} u_i {x_i}^j v_j,
$$
where ${x_i}^j$ stands for the $i$-th component of the $j$-th column of
$X$.  Let $h$ be a continuous function on $I = [$Min $x_i^j$, Max $x_i^j]$.
In order to avoid trivialities, we will assume that $I$ has a nonempty
interior.  For $b=(b_1,\ldots ,b_s)\in\R^s_>$ and $d=(d_0,\ldots ,d_n)
\in\R^{n+1}_>$, let [4, p. 421]
$$
\Hh (b;X;d) = \int_{E_n} \int_{E_{s-1}} h(u\cdot Xv)\phi_b (u)\phi_d (v)\ts
dudv,
$$
$du =du_2 \cdots du_s$, $dv = dv_1 \cdots dv_n$.  Here $\phi_b$ and $\phi_d$
are the Dirichlet densities on $E_{s-1}$ and $E_n$, respectively (see (2.1)).
It is known that for $b\in\R^s_>$,
$$
\Hh (b;X;d) = \int_{E_{s-1}} H(d;u\cdot X) \phi_b(u)\ts du
                                                            \tag 6.1
$$
(see [4, (2.8)]).  In (6.1) $H$ stands for the single Dirichlet average of
$h$ (see (3.1)), $u\cdot X = \{u\cdot x^0, \ldots , u\cdot x^n\}$, $x^0,
\ldots ,x^n$ - the columns of $X$.

     We are in a position to state and prove the following:

\demo{THEOREM 6.1}  Let $d\in\R^{n+1}_>$ and let the vector $b\in\R^s$ be
such that $c\neq 0,-1,\ldots , (c=b_1+\cdots +b_s)$.  If $vol_s([X]) > 0$,
then
$$
\Hh (b;X;d) = \int_{[X]} M(x|d;X) H(b;x)\ts dx,
                                                \tag 6.2
$$
$x=(x_1,\ldots ,x_s)$, $dx = dx_1\cdots dx_s$.
\enddemo

\demo{Proof}  In order to establish (6.2) assume for the moment that
$b\in\R^s_>$.  Application of (3.2) and (3.10) to (6.1) gives
$$
\Hh (b;X;d) = \int_{E_{s-1}} [\int_{[X]}h(u\cdot x)M(x|d;X)\ts dx]
                  \phi_b(u)\ts du.
$$

Interchanging the order of integration and next using (3.1), we obtain the
assertion provided $b\in\R^s_>$.  This restriction can be dropped because
the average $H$ can be continued analytically in the $b$-parameters provided
that $c\neq 0,-1,\ldots$, (see [6,Thm. 6.3--7]).  This completes the proof.
\qd
\enddemo

     Before we state a corollary of Theorem 6.1 let us introduce more
notation.  For $h(z)=z^{-a}$, $(a\in\R)$, the double Dirichlet average of $h$
will be denoted by $\Rr_{-a}$ (cf. [4]).

\demo{COROLLARY 6.2}  ([7])  Let $d\in\R^{n+1}_>$, $b\in\R^s$, and let the
matrix $X$ be such that $0_s \notin [X]$.  Then
$$
   m_{-b}(d;X) = \Rr_{-c}(b;X;d),
                                    \tag 6.3
$$
where $m_{-b}(d;X)$ stands for the moment of order $-c$ of the Dirichlet
spline $M(\cdot |d;X)$.
\enddemo

\demo{Proof} Apply [6,(6.6-5)]
$$ R_{-c}(b;X) = \prod^s_{i=1} x_i^{-b_i}
                                                \tag 6.4
$$
to (6.2) with $h(t) = t^{-c}$. \qd
\enddemo

     Hereafter, we will deal with the hypergeometric functions and polynomials
of several variables.  Appell's hypergeometric function $F_4$ is defined by the
double power series [6, Ex. 6.3-5]
$$
F_4(\a ,\b ; \ga ,\delta ; x_1,x_2) = \sl^{\infty}_{i=0}\sl^{\infty}_{j=0}
      \frac{(\a ,i+j)(\b ,i+j)}{(\ga ,i)(\delta ,i)i!j!}\ts {x_1}^i {x_2}^j,
$$
$\a ,\b ,\ga ,\delta\in\R$, $\ga ,\delta\neq 0,-1,\ldots$,
$|x_1|^{\frac{1}{2}} +|x_2|^{\frac{1}{2}} < 1$.  The following integral
formula [5,p.963]
\newline
$F_4(\a ,\b ; \ga ,\delta ; x_1(1-x_2),x_2(1-x_1)) = $
$$
     \int_0^1 R_{-\a}(d_0,d_1,d_2; u\cdot x^0,u\cdot x^1, u\cdot x^2)
     \phi_b(u)\ts du
                        \tag 6.5
$$
provides the analytic continuation of the $F_4$--series to the region
$\Lambda$ defined by
$$
\Lambda = \{ (x_1,x_2)\in\R^2 : x_1 < 1, x_2 < 1, x_1 + x_2 < 1\}.
$$
In (6.5), $b=(\b ,\ga - \b )$, $d_0 = \ga +\delta -\a -1$, $d_1 = \a + \b - \ga
 -\delta +1$, $d_2 = \delta - \b$, $\phi_b(u)$ is the Dirichlet density on
$E_1$, and $x^0$, $x^1$, and $x^2$ are the columns of $X$, where
$$
X = \bmatrix (1-x_1)(1-x_2) &1-x_1-x_2 &1-x_1 \\
                     1-x_2  &1-x_2     &1 \endbmatrix
                                                      \tag 6.6
$$

\demo{COROLLARY 6.3}  Let $d=(d_0,d_1,d_2)\in\R^3_>$ and let $b=(\b,\ga - \b)
\in\R^2$.  If $vol_2([X]) > 0$, then
$$
F_4(\a ,\b ; \ga ,\delta ; x_1(1-x_2),x_2(1-x_1)) =
         \int_{[X]} M(y|d;X)R_{-\a}(b;y)\ts dy
                                                      \tag 6.7
$$
($y=(y_1,y_2), dy=dy_1dy_2$).  Here, $R_{-\a}$ is the single Dirichlet average
of $h(z)=z^{-\a}$ and the matrix $X$ is given in (6.6).
\enddemo

\demo{Proof}  Apply (6.1) to (6.5) and next use (6.2). \qd
\enddemo

A special case of (6.7) is
$$
F_4(\a ,\b ; \a ,\delta ; x_1(1-x_2),x_2(1-x_1)) = m_{-b}(d;X),
$$
where now $b=(\b ,\a -\b )$ and $d=(\delta -1,\b -\delta +1,\delta - \b)$.
This follows immediately from (6.7) and (6.4).

     We will now deal with Lauricella's $F_B$ function and Lauricella
polynomials.  Let $\a = (\a_1,\ldots ,\a_n)\in\R^n$,
$\b = (\b_1,\ldots ,\b_n)\in\R^n$, $\ga\in\R$\newline
$(\ga \neq 0,-1,\ldots )$, and let $x=(x_1,\ldots x_n)\in\R^n$,
with $|x_i|<1$, for all $i$.  Following [21] we define
$$
F_B(\a ,\b ;\ga ;x) = \sum \frac{(\a ,k)(\b ,k)}{(\ga ,|k|)k!}\ts x^k,
                                                            \tag 6.8
$$
where the summation extends over all multi-indices $k=(k_1,\ldots k_n)\in
\Z^n_+$.  In (6.8) we employ multi-index notation introduced in Section 2.
Also,
$$
   (\a ,k) = \prod^n_{i=1} (\a_i ,k_i).
$$
$(\b ,k)$ is defined in an analogous manner.  When $n=1$, $F_B$ becomes
Gauss' ${}_2F_1$  function.

\demo{COROLLARY 6.4}  Let $d=(\b ,\ga - |\b |)\in\R^{n+1}_>$ and let
$$
X = \bmatrix 1-x_1  &1      &\cdots  &1      &1     \\
                 1  &1-x_2  &\cdots  &1      &1     \\
             \vdots &\vdots &\ddots  &\vdots &\vdots \\
                 1  &    1  &\ldots  &1-x_n  & 1
     \endbmatrix
                                                      \tag 6.9
$$
Then
$$
F_B(\a ,\b ;\ga ;x) = m_{-\a}(d;X)
                                             \tag 6.10
$$
provided that $x_i < 1$ for all $i$.
\enddemo

\demo{Proof}  In the stated domain the entries of $X$ are positive.  Thus
$0_n \notin [X]$.  In order to establish (6.10) we utilize [4,(5.11)] to
obtain
$$
F_B(\a ;\b ;\ga ;x) = \Rr_{-c}(\a ;X;d),
$$
where now $c=\a_1 + \cdots +\a_n$.  This in conjunction with (6.3) gives
the assertion. \qd
\enddemo

     Lauricella polynomials $L_j(x)$, $(j\in\Z^n_+, x\in\R^n)$ are defined
in the following way:
$$
L_j(x) = F_B(-j,\b ;\ga ;x).
$$
These polynomials play an important role in the study of coherent states
(cf. [16]).  On account of (6.10)
$$
L_j(x) = m_j(d;X)
                     \tag 6.11
$$
where the vector $d$ and the matrix $X$ are the same as in Corollary 6.4.

     Two generating functions for the polynomials under discussion can be
derived from (5.9) and (5.11).  Let $\lm = (\lm_1,\ldots ,\lm_n)$,
$e=(1,\ldots ,1)$, \ $\lm ,e\in\R^n$, and let
$d=(\b ,\ga - |\b |)\in\R^{n+1}$.  Then
$$
exp(\lm\cdot e)S(d; -\lm_1x_1,\ldots ,-\lm_nx_n,0) =
                                          \sum \frac{\lm^j}{j!} L_j(x)
                                                               \tag 6.12
$$
$j\in\Z^n_+$.  If
$$
\text{Max}\{|\lm\cdot e - \lm_1x_1|,\ldots ,|\lm\cdot e - \lm_nx_n|,
|\lm\cdot e|\} < 1,
$$
then
$$
    R_{-a}(d;Y) = \sum \lm^j\frac{(a,|j|)}{j!} L_j(x), \tag 6.13
$$
$a\in\R$, $j\in\Z^n_+$, $Y=\{1-\lm\cdot e + \lm_1x_1, \ldots ,
1-\lm\cdot e + \lm_nx_n,1-\lm\cdot e \}$.

For the proof of (6.12) we replace $b$ by $d$ in (5.9) and next use (6.11),
(6.9), and [6,(5.8-3)].  This gives the desired result provided that
$d\in\R^{n+1}_>$.  The latter restriction can be dropped because the $S$--
function can be continued analytically in the $d$-parameters ([6,Corollary
6.3-3]).  The generating function (6.13) can be derived from (5.11) by the
same means.  Feinsilver's generating function [15] can be obtained from
(6.13) by letting $a=\ga$ and then using (4.5).  Recall that $R_0 = 1$.  It is
not hard to show that the $R$-hypergeometric function $R_{-a}$ in (6.13) is a
multiple of the Lauricella function of the fourth kind.  We have

$$
R_{-a}(d;Y) =(1-\lm\cdot e)^{-a} F_D(a,\b ; \ga ; z_1,\ldots ,z_n),
$$
where
$$
z_i = \frac{-\lm_ix_i}{1-\lm\cdot e},
$$
$i=1,2,\ldots ,n$.  We omit further details.

     Our next goal is to establish a recurrence formula obeyed by \newline
Lauricella polynomials
\newline
$ (\ga + |k|)L_{k+d_m}(x) - [\ga (1-w_mx_m) + |k|]L_k(x)$
$$
                              +\sl^n_{l=1} k_l\e_{lm}[L_{k-d_l}(x)
                                 - L_{k-d_l+d_m}(x)] = 0
                                                            \tag 6.14
$$
$m=1,2, \ldots ,n$.  In (6.14), $\ga\in\R_>$, $k=(k_1,\ldots ,k_n)\in\Z^n_+$,
$d_m$ stands for the $m$th coordinate vector in $\R^n$, similarly $d_l$,
$w_m = \frac{\b_m}{\ga}$, $\b_m > 0$, $1\le m\le n$,
$$
\e_{lm} = \cases     1 & \text{if \ }  l\neq m \\
                 1-x_m & \text{if \ }  l=m. \endcases
                                                               \tag 6.15
$$
Here we adopt the convention that $L_k(x)=0$ if $-k_m\in\N$ for some $m$.  In
order to establish the recursion (6.14) we derive first a recurrence formula
for the moments of multivariate Dirichlet splines with $s=n$ and $X$ given in
(6.9).  We have
\newline
$(c+|\b |)m_{\b +d_m}(b;X) - [c(1-w_mx_m)+|\b |]m_{\b}(b;X) $
$$
  + \sl^n_{l=1} \b_l \e_{lm} [m_{\b-d_l}(b;X) - m_{\b - d_l +d_m}(b;X)] = 0,
                                                            \tag 6.16
$$
$b\in\R^{n+1}$, $c=b_1 + \cdots + b_{n+1}$, $\b\in\R^n$, $|\b |=\b_1 + \cdots
\b_n$.  The recursion (6.14) now follows from (6.16) by letting $b=d$, $\b
= k\in\Z^n_+$, and using (6.11).  To complete the proof we need to establish
(6.16).  To this aim we increase the indices of summation in (5.2) and (5.3) by
one unit.  Next we let $s=n$ and solve the resulting linear system for
$m_{\b}(b+e_m;X)$, $1\le m\le n+1$.  Let us note that the assumption
$vol_n([X]) > 0$ is equivalent to $x_m\neq 0$, $1\le m\le n$.  This assures
uniqueness of the solution.  Subtracting (5.3) from (5.2) we obtain
$$
m_{\b}(b+e_m;X) = \frac{[m_{\b}(b;X) - m_{\b + d_m}(b;X)]}{w_mx_m},
                                                                   \tag 6.17
$$
$1\le m \le n$.  The remaining moment $m_{\b}(b+e_{n+1};X)$ can be found using
(5.2) and (6.17).  To complete the proof of (6.16) we utilize (5.4).  Replacing
the index $j$ by $m$ and next using (6.17), we can easily obtain the assertion.

     We close this section with an inequality for Lauricella polynomials.  To
this end, let $x_i <1$, $1\le i \le n$.  It follows from (6.9) that
$[X] \ub \R^n_>$ in the stated domain.  This in conjunction with (6.11) and
(5.1) provides
$$
L_j(x) = \int_{[X]} y^jM(y|d;X)\ts dy > 0
$$
($j\in\Z^n_+$, $d=(\b ,\ga-|\b |)\in\R^{n+1}_>$, $y=(y_1,\ldots ,y_n)$,
$dy=dy_1,\ldots dy_n$).  A standard argument applied to the last formula gives
$$
[L_j(x)]^2 \le L_{j-k}(x)L_{j+k}(x),
$$
where the vector $k\in\Z_+^n$ is such that $j-k\in\Z_+^n$.  In particular,
if \newline\noindent $k=e_m$ -- the $m$th coordinate vector in $\R^n$, then
$$
[L_j(x)]^2 \le L_{j-e_m}(x)L_{j+e_m}(x),
$$
provided $j-e_m\in\Z_+^n$.  Thus the function $g:\Z_+^n\to\R$, where
$$
 g(j_1,\ldots ,j_n) = L_{j_1,\ldots ,j_n}
$$
is log-convex in each variable separately.

\demo{ACKNOWLEDGEMENT}  The authors wish to thank a referee for
insightful comments made on the first draft of this paper, and the  
suggestions for improving Proposition 5.2.
\enddemo
\newpage

\bigskip
\centerline{REFERENCES}
\bigskip
\baselineskip=16truept
\parskip=5pt
\item{[1]}  C. de Boor, Splines as linear combinations of
            $B$-splines:  A survey, in: G.G. Lorentz, C.K. Chui,
            and L.L. Schumaker, Eds., Approximation Theory
            II, (Academic Press, New York, 1976), 1--47.

\item{[2]}  C. de Boor, B-form basics, in:  G. Farin, Ed., Geometric
            Modeling, (SIAM, Philadelphia, 1987), 131--148.

\item{[3]}  C. de Boor and K. H\"ollig, Recurrence relations for
            multivariate $B$-splines, Proc. Amer. Math. Soc.
            85(1982), 397--400.

\item{[4]}  B.C. Carlson, Appell functions and multiple averages,
            SIAM J. Math. Anal. 2(1971), 420--430.

\item{[5]}  B.C. Carlson, Appell's function $F_4$ as a double average,
            SIAM J. Math. Anal. 6(1975), 960-965.

\item{[6]}  B.C. Carlson, Special Functions of Applied
            Mathematics, (Academic Press, New York, 1977).

\item{[7]}  B.C. Carlson, B-splines, hypergeometric functions, and Dirichlet
            averages, J. Approx. Theory 67(1991), 311-325.

\item{[8]}  E. Cohen, T. Lyche, and R.F. Riesenfeld, Cones and
            recurrence relations for simplex splines,
            Constr. Approx. 3(1987), 131--141.

\item{[9]}  H.B. Curry and I.J. Schoenberg, On P\'olya frequency
            functions IV.  The fundamental spline functions and
            their limits, J. d'Analyse Math. 17(1967),
            71--107.

\item{[10]}  W. Dahmen, Multivariate $B$-splines -- Recurrence
            relations and linear combinations of truncated
            powers, in:  W. Schempp and K. Zeller, Eds.,
            Multivariate Approximation Theory, (Basel,
            Birkh\"auser, 1979), 64--82.

\item{[11]}  W. Dahmen, On multivariate $B$-splines, SIAM J.
            Numer. Anal. 17(1980), 993--1012.

\item{[12]}  W. Dahmen, Bernstein-B\' ezier representation of polynomial
            surfaces, in:  Extension of B-spline curve algorithm to
            surfaces, Siggraph 86, organized by C. de Boor, (Dallas 1986).

\item{[13]}  W. Dahmen and C.A. Micchelli, Statistical encounters
            with $B$-splines, Contemp. Math. 59(1986),
            17--48.

\item{[14]}  R. Farwig, Multivariate truncated powers and
            $B$-splines with coalescent knots, SIAM J.
            Numer. Anal. 22(1985), 592--603.

\item{[15]}  P. Feinsilver, Heisenberg algebras in the theory of special
            functions, in: Lecture Notes in Physics, vol. 278 (Springer-
            Verlag, Berlin, 1987), 423--425.

\item{[16]}  P. Feinsilver, Orthogonal polynomials and coherent states, in:
            B. Gruber, L.C. Biedernharn, and M.D. Doehner, Eds., Symmetries in
            Science V (Plenum Press, New York, 1991), 159--172.

\item{[17]} H. Hakopian, Multivariate spline functions,
            $B$-spline basis, and polynomial interpolants,
            SIAM J. Numer. Anal. 19(1982), 510--517.

\item{[18]} K. H\"ollig, A remark on multivariate $B$-splines,
            J. Approx. Theory  33(1981), 119--125.

\item{[19]} S. Karlin, C.A. Micchelli, and Y. Rinott,
            Multivariate splines:  A probabilistic perspective,
            J. Multiv. Anal. 20(1986), 69--90.

\item{[20]} P. Kergin, Interpolation of $C^k$ functions, Ph.D.
            thesis, University of Toronto, 1978.

\item{[21]} G. Lauricella, Sulle funzioni ipergeometriche a piu
            variabli, Rend. Circ. Mat. Palermo 7(1893), 111--158.

\item{[22]} C.A. Micchelli, On a numerically efficient method for
            computing multivariate $B$-splines, in:  W. Schempp
            and K. Zeller, Eds., Multivariate Approximation
            Theory, (Basel, Birkh\"auser, 1979), 211--248.

\item{[23]} C.A. Micchelli, A constructive approach to Kergin
            interpolation in $\R^k$: Multivariate $B$-splines and
            Lagrange interpolation, Rocky Mountain J. Math.
            10(1980), 485--497.

\item{[24]} E. Neuman, Moments and the Fourier transform of
            $B$-splines, J. Comp. Appl. Math. 7(1981),
            51--62.

\item{[25]} E. Neuman, On complete symmetric functions, SIAM
            J. Math. Anal. 19(1988), 736--750.

\item{[26]} E. Neuman and J. Pe\v caric$'$, Inequalities
            involving multivariate convex functions,
            J. Math. Anal. Appl. 137(1989), 541--549.

\item{[27]} E. Neuman, Inequalities involving multivariate convex
            functions II,  Proc. Amer. Math. Soc. 109(1990),
            965--974.

\item{[28]} E. Neuman, Dirichlet averages and their applications to
            special functions (in preparation).
\vfill\eject
\item{[29]} M. Sabin, Open questions in the application of
            multivariate $B$-splines, in: T. Lyche and L.L.
            Schumaker, Eds., Mathematical Methods in
            Computer Aided Geometric Design, (Academic Press,
            Boston, 1989), 529--537.

\item{[30]} L.L. Schumaker, Spline Functions: Basic Theory
            (Wiley, New York, 1981).

\item{[31]} G.S. Watson, On the joint distribution of the
            circular serial correlation coefficients,
            Biometrika 4(1956), 161--168.

\newpage
\centerline{LIST OF SYMBOLS}
\smallskip
\centerline{``Moments of Dirichlet Splines and Their}
\centerline{Applications to Hypergeometric Functions''}
\smallskip
\centerline{by}
\smallskip
\centerline{Edward Neuman and Patrick J. Van Fleet}
\baselineskip=18truept
\bigskip

\settabs\+\indent\indent&letters\quad
               &upper case open-face Roman gamma \cr
\+&$\R$   &upper case open-face Roman R \cr
\+&$\ge$  &greater than or equal to \cr
\+&$\in$  &is a member of \cr
\+&$\sum$ &summation sign \cr
\+&$\subset$ &is a subset of \cr
\+&$\N$   &upper case open-face Roman N \cr
\+&$\Z$   &upper case open-face Roman Z \cr
\+&$\b$   &lower case Greek beta \cr
\+&$!$    &factorial \cr
\+&$\alpha$ &lower case Greek alpha \cr
\+&$\le$  &less than or equal to \cr
\+&$>$    &greater than \cr
\+&$\phi$ &lower case Greek phi \cr
\+&$\prod$ &product \cr
\+&$\Gamma$ &upper case Greek gamma \cr
\+&$\int$  &integral sign \cr
\+&$\Cal F$ &upper case script F \cr
\+&$\Cal H$ &upper case script H \cr
\+&$\Cal R$ &upper case script R \cr
\+&$\epsilon$     &lower case Greek epsilon \cr
\+&$\lm$   &lower case Greek lambda \cr
\+&$\gamma$ &lower case Greek gamma \cr
\+&$\delta$ &lower case Greek delta \cr
\+&$\frac{\pt}{\pt x}$  &partial derivative \cr
\+&$\Lambda$ &upper case Greek lambda \cr
\+&$\eta$ &lower case Greek eta \cr

\bye